\documentclass{article}

\usepackage{inputenc}[utf8]
\usepackage[OT1]{fontenc}

\usepackage{amsmath,amssymb}
\usepackage{amsthm}
\usepackage[pdftex,bookmarks=true]{hyperref}
\usepackage[noadjust]{cite}
\usepackage[capitalise, noabbrev, nameinlink]{cleveref}
\usepackage{enumitem}
\usepackage{tikz}
\tikzset{
every node/.style={draw, circle, inner sep=2pt}, 
edgesign/.style={transform shape, scale=0.6, draw=none, circle, inner sep=1pt, midway, fill=white}, 
textonly/.style={draw=none, rectangle}
}

\usepackage{soul}
\usepackage{cancel}

\newtheorem{theorem}{Theorem}[section]
\newtheorem{lemma}[theorem]{Lemma}
\newtheorem{proposition}[theorem]{Proposition}
\newtheorem{corollary}[theorem]{Corollary}

\theoremstyle{definition}
\newtheorem{definition}[theorem]{Definition}

\newtheorem{remark}[theorem]{Remark}
\newtheorem{example}[theorem]{Example}

\newcommand{\trans}{^\top}

\newcommand{\dunion}{\mathbin{\dot\cup}}
\newcommand{\bzero}{\mathbf{0}}
\newcommand{\bone}{\mathbf{1}}

\newcommand{\bb}{\mathbf{b}}

\newcommand{\by}{\mathbf{y}}
\newcommand{\bz}{\mathbf{z}}

\newcommand{\Col}{\operatorname{Col}}

\newcommand{\spec}{\operatorname{spec}}

\newcommand{\mptn}{\mathcal{S}}
\newcommand{\mptnl}{\mathcal{S}_L}
\newcommand{\mdso}{\ddot{\mathcal{S}}}
\renewcommand{\hom}{\operatorname{Hom}}
\newcommand{\het}{\operatorname{Het}}
\newcommand{\ns}{\mathfrak{S}}
\newcommand{\nw}{\mathfrak{W}}

\title{Congruent quotient matrix, nodal domain, and nodal count}
\author{
Jephian C.-H.~Lin
\thanks{Department of Applied Mathematics, National Yang Ming Chiao Tung University, Hsinchu 300094, Taiwan (jephianlin@gmail.com)}
}
\date{\today}

\begin{document}

\maketitle

\begin{abstract}
We introduce the congruent quotient matrix along with its inequalities on inertia and demonstrate that these inequalities are the essence of the discrete nodal domain theorem, the bounds on nodal counts, and Fiedler's characterization on tree eigenvectors.  Using these techniques, we also extend the existing bounds for the number of nodal domains and nodal counts.
\end{abstract}  


\medskip

\noindent{\bf AMS subject classifications:}
%05C38, %%Paths and cycles
05C50, %%Graphs and linear algebra (matrices, eigenvalues, etc.)
% 05C57, %%Games on graphs
%05C75, %%Structural characterization of families of graphs
% 05C83, %%Graph minors
%05C85, %%Graph algorithms
% 15A03, %%Vector spaces, linear dependence, rank
15A18, %%Eigenvalues, singular values, and eigenvectors
% 15A29, %%Inverse problems
%15B35, %%Sign pattern matrices
15B57, %%Hermitian, skew-Hermitian, and related matrices
%58C15, %%Implicit function theorems; global Newton methods
% 65F18, %%Inverse eigenvalue problems
%68R10, %%Graph theory (including graph drawing)
39A12. %%Discrete version of topics in analysis

\section{Introduction}

Let $A$ be a real symmetric matrix with its eigenvalues $\lambda_1 \leq \cdots \leq \lambda_n$.  Suppose $\by$ is an eigenvector of $A$ with respect to some eigenvalue $\lambda$.  Given the signs of the entries in $A$ and $\by$, how can we estimate the $k$ for $\lambda = \lambda_k$?  Let $n_+(M)$, $n_-(M)$, and $n_0(M)$ be the number of positive, negative, and zero eigenvalues of the real symmetric matrix $M$, respectively.  The triple $(n_+(M), n_-(M), n_0(M))$ is called the \emph{inertia} of $M$.  More precisely, what are the bounds for $n_+(A - \lambda I)$, $n_-(A - \lambda I)$, and $n_0(A - \lambda I)$?  There are many fruitful results on this problem, including the discrete nodal domain theorems \cite{Friedman93, DGLS01, BLS07, Mohammadian16, GL23}, the bounds on nodal counts \cite{Berkolaiko13, Berkolaiko08, BRS12, XY12, CdV13, AG23, AU24, AMU25, BBG25}, and Fiedler's pioneer work on tree eigenvectors \cite{Fiedler75a}.  In this paper, we introduce the notion of congruent quotient matrices to unify the proofs of many theorems and extend the bounds.

We consider finite and undirected simple graphs.  Let $\mathbb{R}^{X\times Y}$ be the set of real matrices whose rows and columns are indexed by $X$ and $Y$, respectively.  Similarly, $\mathbb{R}^Y$ is the set of vectors whose entries are indexed by $Y$.  The $i,j$-entry of a matrix $A$ is denoted by $(A)_{i,j}$, and the $i$-th entry of a vector $\by$ is denoted by $(\by)_i$.

Let $G = (V,E)$ be a graph.  The set $\mptn(G)$ consists of all real symmetric matrices in $\mathbb{R}^{V\times V}$ whose $i,j$-entry, $i\neq j$, is nonzero whenever $\{i,j\}\in E(G)$; the diagonal entries can be any real numbers.  Define 
\[
    \begin{aligned}
    \mdso(G) &= \{A\in\mptn(G): \text{off-diagonal entries} \leq 0\} \text{ and}\\
    \mptnl(G) &= \{A\in\mdso(G): A\bone = \bzero\}.
    \end{aligned}
\]
Here $\bone$ and $\bzero$ are the all-ones vector and the zero vector in the corresponding space.  We call a matrix in $\mdso(G)$ a \emph{discrete Schr\"odinger operator} of $G$ and a matrix in $\mptnl(G)$ a \emph{weighted Laplacian matrix} of $G$.

Let $G = (V,E)$ be a graph and $A\in\mdso(G)$.  Let $\by\in\mathbb{R}^V$ be an eigenvector of $A$ with respect to $\lambda$.  Let $V_+$, $V_-$, $V_0$ be the sets of vertices $v\in V$ such that $(\by)_v$ is positive, negative, and zero, respectively.  A \emph{strong nodal domain} is the vertex set of a component in the induced subgraphs $G[V_+]$ or $G[V_-]$, and the total number of strong nodal domains is denoted by $\ns(\by)$.  In contrast, a \emph{weak nodal domain} is the vertex set of a component in the induced subgraphs $G[V_+\dunion V_0]$ or $G[V_-\dunion V_0]$, and its number is denoted by $\nw(\by)$.  The discrete nodal domain theorem \cite{DGLS01} states that 
\begin{equation}
\label{eq:nd}
    \ns(\by) \leq n_-(A - \lambda I) + n_0(A - \lambda I) \text{ and }
    \nw(\by) \leq n_-(A - \lambda I) + 1,
\end{equation}
where the latter inequality requires $G$ being connected.
The discrete nodal domain theorem is an analog of the Courant's nodal domain theorem \cite{DGLS01}.  Berkolaiko \cite{Berkolaiko08} provided a lower bound for $\ns(\by)$ when $\by$ is nonvanishing, and Xu and Yau \cite{XY12} extended this bound to general $\by$.  Mohammadian \cite{Mohammadian16} introduced a new version of the discrete nodal domain theorem for general matrices $A\in\mptn(G)$.  However, its bound for $\nw(\by)$ can be weaker than the original bound when $A\in\mdso(G)$.  Recently, Ge and Liu \cite{GL23} fixed this issue by redefining the weak nodal domain.

On the other hands, counting the number of edges is also of interest.  Under the same setting of $G$, $A$, and $\by$, define 
\[
    \begin{aligned}
        \hom(A,\by) &= \{\{u,v\}\in E(G): (\by)_u(A)_{u,v}(\by)_v < 0\},  \text{ and} \\
        \het(A,\by) &= \{\{u,v\}\in E(G): (\by)_u(A)_{u,v}(\by)_v > 0\}.
    \end{aligned}
\]
The edges in $\hom(A,\by)$ are said to be \emph{homogeneous}, while the edges in $\het(A,\by)$ are \emph{heterogeneous}.  By definition, $\het(A,\by) = \hom(-A,\by)$.  The number of heterogeneous edges, $|\het(A,\by)|$, is known as \emph{nodal count}.  When $\lambda$ is a simple eigenvalue and $\by$ is nonvanishing, it was shown that 
\begin{equation}
\label{eq:nc}
    n_-(A - \lambda I) \leq |\het(A,\by)| \leq n_-(A - \lambda I) + \beta(G); 
\end{equation}
Here, for a graph $G$ on $n$ vertices with $c$ components, $\beta(G) = |E(G)| - n + c$ is the \emph{cyclomatic number} of $G$, which is the number of extra edges of $G$ comparing to its spanning forests.  The result was first established for $A\in\mdso(G)$ and then proved to be true also for $A\in\mptn(G)$; see \cite{Berkolaiko13, Berkolaiko08, BRS12, CdV13}.  Xu and Yau \cite{XY12} focused on $A\in\mdso(G)$ and extended the two bounds by relaxing the simplicity of $\lambda$ and the nonvanishing condition on $\by$.

Back in 1975, Fiedler provided an elegant connection between the inertia of $A - \lambda I$ and $\by$ when $A\in\mptn(G)$ and $G$ is a tree \cite{Fiedler75a}.  When $\by$ is nonvanishing, Fiedler shows that $n_-(A - \lambda I) = |\het(A,\by)|$ and $n_0(A - \lambda I) = 1$.  Note $\ns(\by) = |\het(A,\by)| + 1$ when $G$ is a tree, so $\ns(\by) = n_-(A - \lambda I) + n_0(A - \lambda I)$.  In conclusion, when $G$ is a tree and $\by$ is nonvanishing, Fiedler not only shows Inequalities \eqref{eq:nd} and \eqref{eq:nc} are tight, but also proves that $\lambda$ is a simple eigenvalue, whereas \cref{eq:nc} takes this as an assumption.  When $\by$ contains some zero entries, Fiedler also shows
\[
    |\het(A,\by)| + p \leq n_-(A - \lambda I) \text{ and } |\hom(A,\by)| + p \leq n_+(A - \lambda I),
\]
where $p$ is the number of vertices $v$ such that $(\by)_v = 0$ and $(\by)_u \neq 0$ for some neighbor of $u$.  

These results are scattered in the literature and provides various aspects to the main problem of estimating the inertia by the signs of $A$ and $\by$.  Many beautiful proofs are given, but it seems that there is a lack of a unified framework to explain these results.  In this paper, we introduce the congruent quotient matrix along with its inequalities and show that these inequalities are the essence of all the abovementioned inequalities, leading to improvements to the bounds.

The paper is organized as follows.  \cref{sec:congruent} introduces the congruent quotient matrix and its inequalities.  They are fundamental linear algebra facts, but we will see that the congruent quotient matrix is the essence of the proofs of all results in this paper.  Other necessary tools are also introduced in this section.  In \cref{sec:strong}, we prove the strong nodal domain theorem and the bounds for nodal counts and extend the existing bounds, which lead to Fiedler's result on nonvanishing tree eigenvectors.  In \cref{sec:parter}, we utilize the Parter vertices to further upgrade the bounds, leading to Fiedler's full result on tree eigenvectors.  We then define the weak nodal domains in \cref{sec:weak}.  Our definition covers the weak nodal domains defined in \cite{Mohammadian16, GL23} as special cases.  Finally, we conclude the paper in \cref{sec:conclusion}.

Lastly, we introduce some notation and terminologies.  Let $A$ be a matrix in $\mathbb{R}^{X\times Y}$ and $\by\in\mathbb{R}^Y$.  For $\alpha\subseteq X$ and $\gamma\subseteq Y$, $A[\alpha,\gamma]$ stands for the submatrix induced on rows in $\alpha$ and columns in $\gamma$, and $\by[\gamma]$ stands for the subvector induced on entries in $\gamma$.  We let $A[\alpha] = A[\alpha,\alpha]$ for short.  Let $G$ be a graph and $e = \{u,v\} \in E(G)$.  If we \emph{contract} the edge $e$, it means that we identify the vertices $u$ and $v$ into one vertex and remove any loop or multiedges.  A \emph{minor} of a graph $G$ is a graph obtained from $G$ by a sequence of edge deletion, vertex deletion, and contraction.  The neighborhood of a vertex $v$ in $G$ is $N_G(v) = \{u: \{u,v\}\in E(G)\}$.

\section{Congruent quotient matrix and necessary tools}
\label{sec:congruent}

We recall that two real symmetric matrices $A$ and $B$ are \emph{congruent} if $B = Q\trans AQ$ for some invertible matrix $Q$, and it is known that congruent matrices have the same inertia; see, e.g., \cite{Zhang26}.  We note that $Q$ being invertible is important, even though the formula $Q\trans AQ$ does not have $Q^{-1}$ involved.  For example, one may choose $Q = O$ and observe that $A$ and $Q\trans AQ$ have different inertia for any nonzero $A$.

\begin{lemma}
\label{lem:congruent}
Let $M$ be an $n\times n$ real symmetric matrix and $S$ an $n\times s$ matrix whose columns are linearly independent.  Let $B = S\trans MS$.  Then the following hold.
\begin{enumerate}[label={\rm(\arabic*)}]
\item $n_-(B) \leq n_-(M)$.
\item $n_-(B) + n_0(B) \leq n_-(M) + n_0(M)$.
\end{enumerate}
\end{lemma}
\begin{proof}
Since the columns of $S$ are linearly independent, one may find an invertible matrix $Q$ such that the first $s$ columns of $Q$ form $S$.  Thus, $M$ and $Q\trans MQ$ are congruent and have the same inertia.  Next, observe that $B = S\trans MS$ is a principal submatrix of $Q\trans MQ$.  By the Cauchy interlacing theorem (see, e.g., \cite{Bapat14}), we have
\begin{enumerate}[label={\rm(\arabic*)}]
\item $n_-(B) \leq n_-(Q\trans MQ) = n_-(M)$, and
\item $n_-(B) + n_0(B) \leq n_-(Q\trans MQ) + n_0(Q\trans MQ) = n_-(M) + n_0(M)$.
\end{enumerate}
This completes the proof.
\end{proof}

We call $B = S\trans MS$ the \emph{congruent quotient matrix} of $M$ with respect to $S$ and demonstrate that \cref{lem:congruent} is the essence of the proofs for theorems in \cref{sec:strong,sec:weak}.

Next, we define the strong nodal domain, which is the same as that in \cite{Mohammadian16, GL23}

\begin{definition}
\label{def:domain}
Let $G$ be a graph, $A\in\mptn(G)$, and $\by\in\mathbb{R}^{V(G)}$.  Define
\[
    \begin{aligned}
        V_{\by= 0} &= \{v\in V(G): (\by)_v = 0\}, \text{ and}\\
        V_{\by\neq 0} &= \{v\in V(G): (\by)_v \neq 0\}.
    \end{aligned}
\]
When $\by$ is clear from the context, we write $G_0 = G[V_{\by = 0}]$ and $G_* = G[V_{\by\neq 0}]$ as the induced subgraphs of $G$.  

A \emph{strong (nodal) domain} is the vertex set of a component of the graph $(V_{\by\neq 0},
\hom(A,\by))$.  A \emph{zero domain} is the vertex set of a component of the induced subgraph $G_0$.  
The number of strong domains is denoted by $\ns(A,\by)$.  When $A$ is clear from the context, we simply write $\ns(\by)$ for $\ns(A,\by)$.
\end{definition}

\begin{figure}[h]
\begin{center}
\begin{tikzpicture}
\node[label={above:$1$}] (1) at (0,0) {};
\node[label={left:$2$}] (2) at (150:1.5) {$+$};
\node[label={left:$3$}] (3) at (210:1.5) {$+$};
\node[label={right:$4$}] (4) at (30:1.5) {$+$};
\node[label={right:$5$}] (5) at (-30:1.5) {$+$};
\draw (4) --node[edgesign]{$+$} (1) --node[edgesign]{$-$} (3) --node[edgesign]{$+$} (2) --node[edgesign]{$+$} (1) --node[edgesign]{$-$} (5);
\draw[very thick] (4) --node[edgesign]{$-$} (5);
\end{tikzpicture}
\end{center}
\caption{A graph $G$ with three strong domains $\{2\}$, $\{3\}$, and $\{4,5\}$, where the homogeneous edges are highlighted.}
\label{fig:sd}
\end{figure}
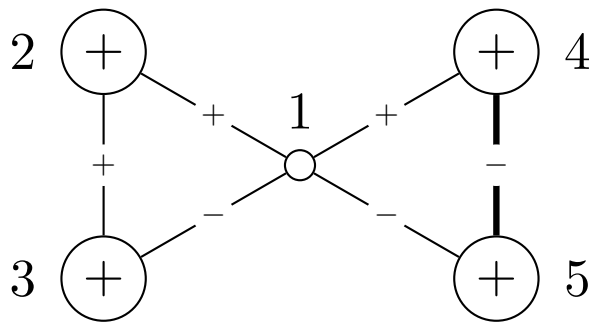

\begin{example}
\label{ex:sd}
Let $G$ be the graph as shown in \cref{fig:sd}.  Let 
\[
    A = \begin{bmatrix}
        0 & 1 & -1 & 1 & -1 \\
        1 & -1 & 1 & 0 & 0 \\
        -1 & 1 & -1 & 0 & 0 \\
        1 & 0 & 0 & 1 & -1 \\
        -1 & 0 & 0 & -1 & 1
    \end{bmatrix}
    \text{ and }
    \by = \begin{bmatrix}
        0 \\ 1 \\ 1 \\ 1 \\ 1
    \end{bmatrix}.
\]
In \cref{fig:sd}, the signs of $A$ are marked on the edges, and the signs of $\by$ are marked inside the vertices.  By definition, $\{4,5\}$ is the only homogeneous edge, and there are three strong domains $\{2\}$, $\{3\}$, and $\{4,5\}$.  
\end{example}

By definition, when $A\in\mdso(G)$ is a discrete Schr\"odinger operator, the strong domain here aligns with the strong nodal domain considered in \cite{DGLS01}.

Let $G$ be a graph on $n$ vertices with $c$ components.  We recall that any spanning forest of $G$ has $n - c$ edges.  

\begin{proposition}
\label{prop:spforest}
Let $G$ be a graph and $A\in\mptn(G)$.  Let $\by\in\mathbb{R}^{V(G)}$.  Then the following holds.
\begin{enumerate}[label={\rm(\arabic*)}]
\item $n - z - \ns(A,\by) \leq |\hom(A,\by)|$.
\item $n - z - \ns(-A,\by) \leq |\het(A,\by)|$.
\item $|\hom(A,\by)| + |\het(A,\by)| = |E(G_*)| = n - z - c_* + \beta_*$.
\end{enumerate}
Here $n = |V(G)|$, $z = |V_{\by= 0}|$, $c_*$ is the number of components of $G_*$, and $\beta_* = \beta(G_*)$ is the cyclomatic number of $G_*$.  
\end{proposition}
\begin{proof}
The graph $(V_{\by\neq 0}, \hom(A,\by))$ has $n - z$ vertices and $\ns(A,\by)$ components, so any of its spanning trees has $n - z - \ns(A,\by)$ edges, which gives the first inequality.  The second inequality is obtained by applying the first inequality to $-A$ and $\by$.  The last equality comes from the fact that $\hom(A,\by) \dunion \het(A,\by)$ is a partition of $E(G_*)$.
\end{proof}

Finally, we recall some classical results regarding the weighted Laplacian matrices.

\begin{proposition}[See, e.g., \cite{Bapat14}]
\label{prop:psd}
Let $G$ be a graph with $c$ components and $A\in\mptnl(G)$.  Then $A$ is a positive semidefinite matrix with $n_0(A) = c$.  
\end{proposition}

\cref{prop:psd} is not true for general matrices.  For example, consider the $4\times 4$ matrix $A = \begin{bmatrix} (-1)^{i + j - 1} \end{bmatrix}$, which is a matrix in $\mptn(K_4)$.  Then $A$ is not a positive semidefinite matrix.  Moreover, $n_0(A) = 3$, which is greater than the number of components, $c = 1$, of $K_4$.  
% However, with the condition $A\bone = \bzero$, we have an inequality.

% \begin{proposition}
% \label{prop:gpsd}
% Let $G$ be a graph with $c$ components.  If $A\in\mptn(G)$ and $A\bone = \bzero$, then $n_0(A) \geq c$.
% \end{proposition}
% \begin{proof}
% Let $X_1, \ldots, X_c$ be the vertex sets of the components of $G$.  Then $A$ is the direct sum of $A[X_1], \ldots, A[X_c]$.  Since $A\bone = \bzero$, $\ker(A[X_k])$ contains the all-ones vector in $\mathbb{R}^{X_k}$ for each $k = 1, \ldots, c$, which means $n_0(A[X_k]) \geq 1$.  Therefore, $n_0(A) \geq c$.
% \end{proof}

\section{Strong nodal domain and nodal count}
\label{sec:strong}

In this section, we provide upper and lower bounds for the number of strong nodal domains, $\ns(A,\by)$, and the number of heterogeneous edges, $|\het(A,\by)|$.

Recall that the nodal domain theorem \cite{DGLS01} states that $\ns(A,\by) \leq n_-(A - \lambda I) + n_0(A - \lambda I)$ when $A\in\mdso(G)$, which was then proved to be true also for general matrices $A\in\mptn(G)$ in \cite{Mohammadian16, GL23}.  On the other hand, Berkolaiko shows a lower bound of $\ns(A,\by)$ for $A\in\mdso(G)$ and nonvanishing $\by$ \cite{Berkolaiko08}.  Later, Xu and Yau extended the lower bound to 
\[
    n_-(A - \lambda I) + n_0(A - \lambda I) - z - \beta(G) \leq \ns(A,\by),
\]
where $z = |V_{\by = 0}|$ and $\by$ may contain zero entries \cite{XY12}.

\cref{thm:jsndt} not only enhance the bounds on both sides but also generalize Xu and Yau's lower bound for discrete Schr\"odinger operators to general matrices $A\in\mptn(G)$.  Meanwhile, it demonstrates how the congruent quotient matrix streamlines the proof.

\begin{theorem}
\label{thm:jsndt}
Let $G$ be a graph and $A\in\mptn(G)$.  If $\by$ is an eigenvector of $A$ with respect to $\lambda$, then  
\[
    n_- + \max\{c_*, n_0\} - z - \beta_* \leq \ns(A,\by) \leq n_- + \min\{c_*, n_0\}.
\]
Here $n_- = n_-(A - \lambda I)$, $n_0 = n_0(A - \lambda I)$, $z = |V_{\by = 0}|$, $c_*$ is the number of components of $G_*$, and $\beta_* = \beta(G_*)$ is the cyclomatic number of $G_*$.  
\end{theorem}
\begin{proof}
We first prove $\ns(A,\by) \leq n_- + \min\{c,n_0\}$.  Let $A_\lambda = A - \lambda I$.  Let $X_1, \ldots, X_q$ be the strong domains with respect to $A$ and $\by$, where $q = \ns(A,\by)$.  Define $\by_i$ as the vector obtained from $\by$ by keeping the same value for entries in $X_i$ and setting it to zero for others.  Thus, $\{\by_1, \ldots, \by_q\}$ is a linearly independent set since they are nonzero vectors with mutually disjoint supports.  Now let $S$ be the $n\times q$ matrix whose $i$-th column is $\by_i$ and consider $B = S\trans A_\lambda S$.  Note that the entry of $B$ has the form $\by_i\trans A_\lambda\by_j$, which can be written as the sum of $(\by)_u (A)_{u,v} (\by)_v$ over edges $\{u,v\}$ between $X_i$ and $X_j$.  Since edges connecting different strong domains have to be heterogeneous, that is, $(\by)_u A_{u,v} (\by)_v \geq 0$, we know that the off-diagonal entries of $B$ is zero or positive.  More precisely, for $i\neq j$, we have $(B)_{i,j} > 0$ if there is an edge connecting $X_i$ and $X_j$, and $(B)_{i,j} = 0$ otherwise.  Moreover, 
\[
    B\bone = S\trans A_\lambda S\bone = S\trans A_\lambda \by = \bzero.
\]
In summary, $-B$ is a matrix in $\mptn_L(H)$, where $H$ be the graph obtained from $G_*$ by contracting all homogeneous edges.  Note that $H$ has $c_*$ components, since contraction does not change the number of components.

By \cref{prop:psd}, $B$ is a $q\times q$ negative semidefinite matrix with $n_0(B) = c_*$, so we obtain two inequalities
\[
    \begin{aligned}
        \ns(A,\by) - c_* &= n_-(B) \leq n_-, \\
        \ns(A,\by) &= n_-(B) + n_0(B) \leq n_- + n_0
    \end{aligned}
\]
by \cref{lem:congruent}, which can be combined as $\ns(A,\by) \leq n_- + \min\{c_*,n_0\}$.

Next we show that the lower bound can be obtained from the upper bound together with \cref{prop:spforest}.  The three items in \cref{prop:spforest} imply 
\[
    n - z - \ns(A,\by) + n - z - \ns(-A,\by) \leq n - z - c_* + \beta_*,
\]
which leads to 
\[
    n - z - \ns(-A,\by) + c_* - \beta_* \leq \ns(A,\by).
\]
By applying \cref{thm:jsndt} to $-A$ and $\by$ with respect to the eigenvalue $-\lambda$, we have 
\[
    \begin{aligned}
    \ns(-A,\by) &\leq n_-(-A_\lambda) + \min\{c_*, n_0(-A_\lambda)\} \\ 
    &= n_+(A_\lambda) + \min\{c_*, n_0(A_\lambda)\} \\
    &= n_+ + \min\{c_*, n_0\} \\
    &= n_+ + n_0 + c_* - \max\{c_*, n_0\},
    \end{aligned}
\]
where $n_+ = n_+(A_\lambda)$.
Combining these inequalities, we have the desired result.
\end{proof}

% \begin{remark}
% We note that the replacement of $n_0$ with $\max\{c,n_0\}$ or $\min\{c,n_0\}$ tightened the bounds on both sides, and we will see later that there are examples where \cref{thm:jsndt} indeed provides tighter bounds.  Moreover, by definition we have $\beta_* \leq \beta(G)$, which leads to another improvements.
% \end{remark}

The following examples show that $c_*$ and $n_0$ are not comparable.  When $A\in\mdso(G)$, we recall the upper bound $\ns(A,\by) \leq n_- + n_0$ in \cite{DGLS01} and the lower bound $\ns(A,\by) \geq n_- + n_0 - z - \beta(G)$ in \cite{XY12}.  Comparing to these results, we see that \cref{thm:jsndt} enhances the upper bound when $c_* < n_0$ and the lower bound when $c_* > n_0$.

\begin{example}
\label{ex:k3}
Let $G = K_3$.  Let 
\[
    A = \begin{bmatrix}
        -1 & -1 & -1 \\
        -1 & -1 & -1 \\
        -1 & -1 & -1 \\
    \end{bmatrix}
    \text{ and }
    \by = \begin{bmatrix}
        0 \\ 1 \\ -1
    \end{bmatrix}.
\]
Then $\spec(A) = \{-3, 0, 0\}$ and $\by$ is an eigenvector of $A$ with respect to $\lambda = 0$.  Applying \cref{thm:jsndt} to $A$ and $\by$, we have $n_- = 1$, $n_0 = 2$, and $c_* = 1$.  This gives a tight upper bound 
\[
    2 = \ns(A,\by) \leq 1 + 1 = 2.
\]
Note that $c_* < n_0$ in this example.
\end{example}

\begin{example}
\label{ex:k14}
Let $G$ be the graph obtained from $K_{1,4}$ by adding an edge.  Let 
\[
    A = \begin{bmatrix}
        0 & -1 & -1 & -1 & -1 \\
        -1 & -1 & -1 & 0 & 0 \\
        -1 & -1 & -1 & 0 & 0 \\
        -1 & 0 & 0 & 0 & 0 \\
        -1 & 0 & 0 & 0 & 0 \\
    \end{bmatrix}
    \text{ and }
    \by = \begin{bmatrix}
        0 \\ 1 \\ -1 \\ 1 \\ -1
    \end{bmatrix}.
\]
Then $\spec(A)$ is approximately $\{-2.9, -0.8, 0, 0, 1.7\}$, and $\by$ is an eigenvector of $A$ with respect to $\lambda = 0$.  Applying \cref{thm:jsndt} to $A$ and $\by$, we have $n_- = 2$, $n_0 = 2$, $z = 1$, $c_* = 3$, and $\beta_* = 0$.  This gives a tight lower bound
\[
    4  = \ns(A,\by) \geq 2 + 3 - 1 - 0 = 4.
\]
Note that $c_* > n_0$ in this example.
\end{example}

% \begin{remark}
% While $c_*$ and $n_0$ are not comparable in general, we know by \cref{prop:psd,prop:gpsd} that $c_* \leq n_0$ when $z = 0$, and the equality holds when $A\in\mdso(G)$.
% \end{remark}

\cref{thm:f75a1} is by Fiedler.  Here we show that it can be viewed as a consequence of \cref{thm:jsndt}.

\begin{theorem}[Fiedler~\cite{Fiedler75a}]
\label{thm:f75a1}
Let $T$ be a tree and $A\in\mptn(T)$.  If $\by$ is a nonvanishing eigenvector with respect to $\lambda$, then
\[
    \ns(A,\by) = |\het(A,\by)| + 1 = n_- + 1
\]
and $\lambda$ is a simple eigenvalue, where $n_- = n_-(A - \lambda I)$.
\end{theorem}
\begin{proof}
The fact that $T$ is a tree implies $\ns(A,\by) = |\het(A,\by)| + 1$.
We apply \cref{thm:jsndt} to $T$, $A$, and a nonvanishing $\by$.  Since $T$ is a tree and $\by$ is nonvanishing, we have $c_* = 1$, $z = 0$, and $\beta_* = 0$.  This leads to the inequalities 
\[
    n_- + n_0 \leq \ns(A,\by) \leq n_- + 1,
\]
where $n_- = n_-(A - \lambda I)$ and $n_0 = n_0(A - \lambda I)$.  Equivalently, this means $n_0 = 1$ and $\ns(A,\by) = n_- + 1$.
\end{proof}

Next we derive the bounds on nodal count from \cref{thm:jsndt}.

\begin{corollary}
\label{cor:jnc}
Let $G$ be a graph and $A\in\mptn(G)$.  If $\by$ is an eigenvector of $A$ with respect to $\lambda$, then  
% \[
%      n_- + n_0 - z - \min\{c,n_0\} \leq |\het(A,\by)| \leq n_- + n_0 +\beta_* - \max\{c,n_0\}.
% \]
\[
     n_- + \max\{c_*,n_0\} - c_* - z \leq |\het(A,\by)| \leq n_- + \min\{c_*,n_0\} - c_* +\beta_*.
\]
Here $n_- = n_-(A - \lambda I)$, $n_0(A - \lambda I)$, $z = |V_{\by = 0}|$, $c_*$ is the number of components of $G_*$, and $\beta_* = \beta(G_*)$ is the cyclomatic number of $G_*$.  
\end{corollary}
\begin{proof}
This follows from \cref{thm:jsndt} and \cref{prop:spforest}.  Let $A_\lambda = A - \lambda I$.  By \cref{prop:spforest}, we have 
\[
    \begin{aligned}
        |\het(A,\by)| & \geq n - z - \ns(-A_\lambda,\by) \\
        & \geq n - z - n_-(-A_\lambda) - \min\{c_*, n_0(-A_\lambda)\} \\
        &= n - z - n_+(A_\lambda) - \min\{c_*, n_0\} \\
        &= n_- + n_0 - z - \min\{c_*, n_0\} \\
        &= n_- + \max\{c_*,n_0\} - c_* - z,
    \end{aligned}
\]
where the second inequality is by applying \cref{thm:jsndt} to $-A$, $\by$, and $-\lambda$.

On the other hand, by \cref{prop:spforest} we have 
\[
    \begin{aligned}
        |\het(A,\by)| &= n - z - c_* + \beta_* - |\hom(A,\by)| \\
        &\leq n - z - c_* + \beta_* - (n - z - \ns(A,\by)) \\
        &\leq n_- + \min\{c_*,n_0\} - c_* + \beta_*.
        % &= n_- + n_0 + \beta_* - \max\{c,n_0\}
    \end{aligned}
\]
This completes the proof.
\end{proof}

\section{Parter vertex}
\label{sec:parter}

We have seen that \cref{thm:jsndt} implies Fiedler's result \cref{thm:f75a1} when the graph is a tree and the eigenvector is nonvanishing.  In this section, we show that both \cref{thm:jsndt,cor:jnc} may integrate with Parter vertices to produce more refined bounds.

We first state the full version of Fiedler's theorem.

\begin{theorem}[Fiedler~\cite{Fiedler75a}]
\label{thm:f75a2}
Let $T$ be a tree and $A\in\mptn(G)$.  If $\by$ is an eigenvector with respect to $\lambda$, then
\[
    |\het(A,\by)| + p \leq n_- \text{ and } |\hom(A,\by)| + p \leq n_+,
\]
where $n_- = n_-(A - \lambda I)$, $n_+ = n_+(A - \lambda I)$, and $p$ is the number of vertices $v$ such that $(\by)_v = 0$ and $(\by)_u \neq 0$ for some neighbor of $v$.
\end{theorem}

\begin{figure}[h]
\begin{center}
\begin{tikzpicture}
\node[label={below:$1$}] (1) at (0,0) {};
\node[label={below:$2$}] (2) at (-1,0) {};
\node[label={left:$3$}] (3) at (-2,0) {$+$};
\node[label={left:$4$}] (4) at (-2,1) {$-$};
\node[label={above:$5$}] (5) at (0,1) {};
\node[label={left:$6$}] (6) at (-1,2) {$+$};
\node[label={right:$7$}] (7) at (1,2) {$-$};
\node[label={below:$8$}] (8) at (1,0) {};
\node[label={right:$9$}] (9) at (2,0) {$+$};
\node[label={right:$10$}] (10) at (2,1) {$-$};
\draw (3) --node[edgesign]{$-$} (2) --node[edgesign]{$-$} (4);
\draw (6) --node[edgesign]{$-$} (5) --node[edgesign]{$-$} (7);
\draw (9) --node[edgesign]{$-$} (8) --node[edgesign]{$-$} (10);
\draw (2) --node[edgesign]{$-$} (1) --node[edgesign]{$-$} (8);
\draw (1) --node[edgesign]{$-$} (5);
\end{tikzpicture}
\end{center}
\caption{A tree $T$ with three vertices, $v = 2, 5, 8$, such that $(\by)_v = 0$ and $(\by)_u \neq 0$ for some neighbor of $v$.}
\label{fig:fiedler}
\end{figure}
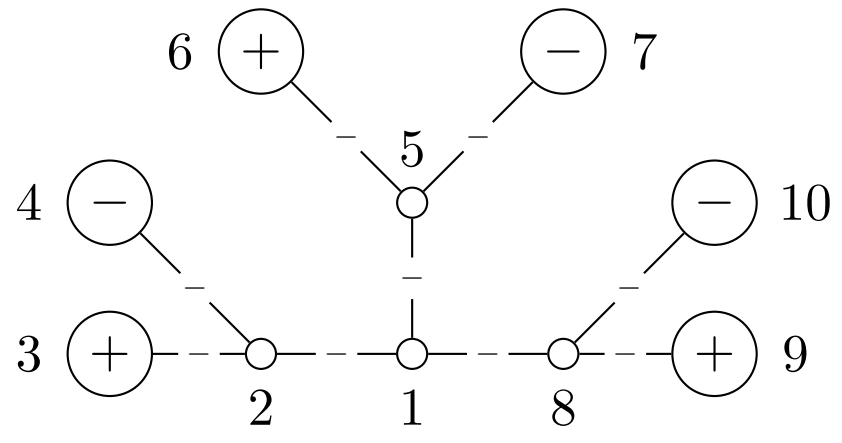

\begin{example}
\label{ex:fiedler}
Let $T$ be the tree as shown in \cref{fig:fiedler}.  Let $A\in\mptn(T)$ be the negative of the adjacency matrix of $T$.  Then 
\[
    \by = (0,0,1,-1,0,1,-1,0,1,-1)\trans
\]
is an eigenvector of $A$ with respect to $0$, and its signs are illustrated in \cref{fig:fiedler} inside the nodes.  A direct computation shows that the inertia of $A$ is $(n_+, n_-, n_0) = (3, 3, 4)$.  Also, there are three vertices, $v = 2, 5, 8$, such that $(\by)_v = 0$ and $(\by)_u \neq 0$ for some neighbor of $v$.  Thus, \cref{thm:f75a2} gives $n_- \geq 3$ since $|\het(A,\by)| = 0$.  Both \cref{thm:jsndt,cor:jnc} give $n_- \geq 6 - \min\{6,n_0\} = 2$, which is not tight even if $n_0$ is known beforehand.  
\end{example}

\begin{definition}
Let $G$ be a graph and $A\in\mptn(G)$.  Let $\by$ be an eigenvector of $A$ with respect to $\lambda$.  Let $v\in V(G)$ and $X\subseteq V(G)$ with $v\notin X$.  Consider the multiset
\begin{equation}
\label{eq:vxsum}
    \{(A)_{v,u}(\by)_u: u\in X,\ \{v,u\}\in E(G)\}
\end{equation}
Define
\begin{itemize}
    \item $s(v,X) = +$ if \cref{eq:vxsum} contains at least one positive value and no negative values,
    \item $s(v,X) = -$ if \cref{eq:vxsum} contains at least one negative value and no positive values,
    \item $s(v,X) = 0$ if \cref{eq:vxsum} is empty or contains only zero values,
    \item $s(v,X) = \circledast$ if \cref{eq:vxsum} contains at least one positive and at least one negative value.
\end{itemize}
When $s(v,X)\in\{+,-\}$, we say $s(v,X)$ is \emph{nonzero}.
\end{definition}

For example, we have $s(2,\{3\}) = -$ and $s(2,\{3,4\}) = \circledast$ in \cref{fig:wdone}.  
The intuition of $s(v,X)$ is straightforward.  If $s(v,X) = +$, then the sum of elements in \cref{eq:vxsum} is always positive.  The other cases are similar.  

\begin{definition}
Let $G$ be a graph and $A\in\mptn(G)$.  Let $\by$ be an eigenvector of $A$ with respect to $\lambda$.  A vertex $v$ is called a \emph{signed zero} if $(\by)_v = 0$ and $s(v,X)$ is nonzero for some component of $G - v$ with vertex set $X$.  
\end{definition}

Let $G$ be a graph and $A\in\mptn(G)$.  Let $\lambda\in\spec(A)$ and $A_\lambda = A - \lambda I$.  Recall that a \emph{Parter vertex} is a vertex $v\in V(G)$ such that $n_0(A_\lambda(v)) = n_0(A_\lambda) + 1$.  Here $M(v)$ stands for the matrix obtained from $M$ by removing the row and the column indexed by $v$.  The Parter vertex is motivated by the Parter--Wiener theorem \cite{Parter60,Wiener84,JLDS03}, and it is an active research area in spectral graph theory \cite{DH22,DdF25,KS08,SP26}.  Note that a Parter vertex is sometimes referred to as an upper vertex; see, e.g., \cite{SNIP,temno,IEPS}.  Now we may justify our definition of signed zeros.

\begin{proposition}
\label{prop:szp}
Let $G$ be a graph and $A\in\mptn(G)$.  Let $\by$ be an eigenvector of $A$ with respect to $\lambda$.  If $v$ is a signed zero, then $v$ is a Parter vertex.  More precisely,
\[
    \begin{aligned}
    n_+(A_\lambda(v)) &= n_+(A_\lambda) - 1, \\ 
    n_-(A_\lambda(v)) &= n_-(A_\lambda) - 1, \text{ and }\\ 
    n_0(A_\lambda(v)) &= n_0(A_\lambda) + 1.
    \end{aligned}
\]
where $A_\lambda = A - \lambda I$.
\end{proposition}
\begin{proof}
Let $v$ be a signed zero and $X_1,\ldots,X_c$ the vertex sets of components of $G - v$.  Without loss of generality, we may assume $v = 1$ and $s(v,X_1)$ is nonzero.  Thus, we write
\[
    A_\lambda = \begin{bmatrix}
        a & \bb_1\trans & \cdots & \bb_c\trans \\
        \bb_1 & A_1 & ~ & ~ \\
        \vdots & ~ & \ddots & ~ \\
        \bb_c & ~ & ~ & A_c
    \end{bmatrix},
\]
where $A_k = A_\lambda[X_k]$ for $k = 1, \ldots, c$.  Since $\by\in\ker(A_\lambda)$ and $(\by)_v = 0$, $\by[X_1]\in\ker(A_1)$.  Since $s(v,X_1)$ is nonzero, $\bb_1\trans\by[X_1] \neq 0$, which implies $\bb_1\notin\Col(A_1)$.  This further implies $\bb\notin\Col(A_\lambda(v))$, where $\bb = (\bb_1\trans, \cdots, \bb_c\trans)\trans$.  As a consequence, $n_0(A_\lambda(v)) = n_0(A_\lambda) + 1$.  By the Cauchy interlacing theorem (see, e.g., \cite{Bapat14}), this also means $n_+(A_\lambda(v)) = n_+(A_\lambda) - 1$ and $n_-(A_\lambda(v)) = n_-(A_\lambda) - 1$.
\end{proof}

Under the setting of Fiedler's result (\cref{thm:f75a2}), a vertex $v$ such that $(\by)_v = 0$ and $(\by)_u \neq 0$ for some neighbor of $v$ is a signed zero since the considered graph is a tree.  In \cite{GL23}, for a graph $G$, a vertex $v$ is called \emph{tree-like} if the number of components of $G - v$ equals the degree of $v$, and it is used to improve the lower bound of $\ns(A,\by)$.  Indeed, if a graph $G$ contains a tree-like vertex $v$ such that $(\by)_v = 0$ and $(\by)_u \neq 0$ for some neighbor of $v$, then $v$ is also a signed zero.  Therefore, the notion of signed zeros generalizes these cases.

\begin{corollary}
\label{cor:szp}
Let $G$ be a graph and $A\in\mptn(G)$.  Let $\by$ be an eigenvector of $A$ with respect to $\lambda$.  If there are $p$ signed zeros, then 
\[
    \begin{aligned}
    n_+(A_\lambda[V']) &= n_+(A_\lambda) - p, \\ 
    n_-(A_\lambda[V']) &= n_-(A_\lambda) - p, \text{ and }\\ 
    n_0(A_\lambda[V']) &= n_0(A_\lambda) + p.
    \end{aligned}
\]
where $A_\lambda = A - \lambda I$ and $V'$ is obtained from $V(G)$ by removing the $p$ signed zeros.
\end{corollary}
\begin{proof}
Let $u$ and $v$ be signed zeros with respect to $G$, $A$, and $\by$.  Let $V_v = V(G) - v$.  The desired results follow inductively from the observation that a signed zero with respect to $G$, $A$, and $\by$ is also a signed zero with respect to $G[V_v]$, $A[V_v]$, and $\by[V_v]$.
\end{proof}

Next we use the Parter vertices to improve the bounds in \cref{sec:strong}.

% \[
%     n_- + \max\{c_*,n_0\} - z - \beta_* \leq \ns(A,\by) \leq n_- + \min\{c_*,n_0\}
% \]
% and 
% \[
%     n_- + \max\{c,n_0\} - c_* - z \leq |\het(A,\by)| \leq n_- + \min\{c_*,n_0\} - c_* + \beta_*.
% \]
\begin{theorem}
\label{thm:jsndtp}
Let $G$ be a graph and $A\in\mptn(G)$.  Let $\by$ be an eigenvector $A$ with respect to $\lambda$.  Then
\[
    n_- + \max\{c_*,n_0 + p\} - z - \beta_* \leq \ns(A,\by) \leq n_- + \min\{c_* - p,n_0\}
\]
and 
\[
    n_- + \max\{c_*,n_0 + p\} - c_* - z \leq |\het(A,\by)| \leq n_- + \min\{c_* - p,n_0\} - c_* + \beta_*.
\]
Here $(n_+,n_-,n_0)$ is the inertia of $A - \lambda I$, $c_*$ is the number of component of $G_*$, $\beta_* = \beta(G_*)$ is the cyclomatic number of $G_*$, and $p$ is the number of signed zeros.
\end{theorem}
\begin{proof}
Let $P$ be the set of signed zeros with respect to $G$, $A$, and $\by$.  Let $V' = V \setminus P$.  We apply \cref{thm:jsndt,cor:jnc} to $G[V']$, $A[V']$, and $\by[V']$.  First, we note that $\ns(G,\by)$, $|\het(A,\by)|$, $c_*$, and $\beta_*$ are parameters related to $G_*$, which is unchanged when the vertices in $P$ are removed, so these parameters stay the same.  It is straightforward to see that $z$ drop by $p$ when the context changes from $G$, $A$, and $\by$ to $G[V']$, $A[V']$, and $\by[V']$.  Finally, by \cref{cor:szp}, we have
\[
    \begin{aligned}
    \ns(A,\by) &= \ns(A[V'],\by[V']) \\
    &\geq (n_- - p) + \max\{c_*, n_0 + p\} - (z - p) - \beta_* \\
    &= n_- + \max\{c_*,n_0 + p\} - z - \beta_*, \\
    \ns(A,\by) &= \ns(A[V'],\by[V']) \\
    &\leq (n_- - p) + \min\{c_*,n_0 + p\} \\
    &= n_- + \min\{c_* - p,n_0\}, \\
    |\het(A,\by)| &= |\het(A[V'],\by[V'])| \\
    &\geq (n_- - p) + \max\{c_*,n_0 + p\} - c_* - (z - p) \\
    &= n_- + \max\{c_*,n_0 + p\} - c_* - z, \\
    |\het(A,\by)| &= |\het(A[V'],\by[V'])| \\
    &\geq (n_- - p) + \min\{c_*,n_0 + p\} - c_* + \beta_* \\
    &= n_- + \min\{c_* - p,n_0\} - c_* + \beta_*.
    \end{aligned}
\]
This completes the proof.
\end{proof}

Applying \cref{thm:jsndtp} to \cref{ex:fiedler}, one may derive bounds $n_- + 3 \leq \ns(A,\by) \leq n_- + 3$.  Thus, the bounds are tight and $n_- = \ns(A,\by) - 3 = 3$.  Indeed, Fiedler's result (\cref{thm:f75a2}) can be viewed as a corollary of \cref{thm:jsndtp}.

\begin{proof}[Proof of \cref{thm:f75a2}]
Recall that a vertex $v$ such that $(\by)_v = 0$ and $(\by)_u \neq 0$ for some neighbor of $v$ is a signed zero.  Let $p$ be the number of such vertices.  Applying \cref{thm:jsndtp} to $T$, $A$, and $\by$, we have 
\[
    \begin{aligned}
        |\het(A,\by)| &\leq n_- + \min\{c_* - p,n_0\} - c_* + \beta_* \\
        &\leq n_- + (c_* - p) - c_* + \beta_* \\
        % &= n_- - \min\{-p, n_0 - c_*\} \\
        &= n_- - p.
    \end{aligned}
\]
since $\beta_* = 0$, which gives the first inequality.  The second inequality comes from applying the first inequality to $-A$ and $\by$.
\end{proof}

\section{Discrete weak nodal domain theorem}
\label{sec:weak}

In this section, we define the weak nodal domain to extend the bound $\nw(\by) \leq n_-(A - \lambda I) + 1$ in the classical discrete nodal domain theorem \cite{DGLS01}.  The congruent quotient matrix is again the key to the proof.  A few more linear algebra tools are needed.

\begin{lemma}
\label{lem:weyl}
Let $M$ be a real symmetric matrix and $P$ a positive semidefinite matrix.  Then the following holds.
\begin{enumerate}[label={\rm(\arabic*)}]
\item $n_-(M + P) \leq n_-(M)$.
\item $n_-(M + P) + n_0(M + P) \leq n_-(M) + n_0(M)$.
\end{enumerate}
\end{lemma}
\begin{proof}
This follows from Weyl's inequality; see, e.g., \cite{Zhang26}.  In this special case, it states that adding a positive semidefinite matrix does not decrease the eigenvalues.
\end{proof}

Recall that a \emph{weighted incidence matrix} of a graph $G$ is a matrix in $\mathbb{R}^{V(G)\times E(G)}$ such that the column corresponding to the edge $\{u,v\}$
is a vector in $\mathbb{R}^{V(G)}$ with exactly two nonzero entries, where one of the $u$-th entry and $v$-th entry is $w_{u,v}$ and the other is $-w_{u,v}$ for some $w_{u,v} > 0$.  Note that if $N$ is a weighted incidence matrix of $G$, then $NDN\trans\in\mptnl(G)$ for any diagonal matrix $D$ with positive diagonal entries.

We note that the union $G_1\cup G_2$ in \cref{lem:dnngamma} stands for the graph $(V(G_1)\cup V(G_2), E(G_1)\cup E(G_2))$, rather than the disjoint union of them.

\begin{lemma}
\label{lem:dnngamma}
Let $G_1$ and $G_2$ be graphs such that $V(G_1) = V(G_2)$ and $G_1\cup G_2$ is connected.  Let $B$ be a matrix of the form 
\[
    B = \begin{bmatrix}
        D & N\trans \\
        N & \Gamma
    \end{bmatrix}
\]
such that $D$ is a diagonal matrix with positive diagonal entries, $N$ is a weighted incidence matrix of $G_2$, and $-\Gamma\in\mptnl(G_1)$.  Then $B$ is a matrix with $n_-(B) = |V(G_1)| - 1$ and $n_0(B) = 1$.
\end{lemma}
\begin{proof}
Observe that $B$ is congruent to 
\[
    \begin{bmatrix}
        I & O \\
        -ND^{-1} & I
    \end{bmatrix}
    \begin{bmatrix}
        D & N\trans \\
        N & \Gamma
    \end{bmatrix}
    \begin{bmatrix}
        I & -D^{-1}N\trans \\
        O & I
    \end{bmatrix}
    =
    \begin{bmatrix}
        D & O \\
        O & \Gamma - ND^{-1}N\trans
    \end{bmatrix}.
\]
Since $-\Gamma\in\mptnl(G_1)$ and $ND^{-1}N\trans\in\mptnl(G_2)$, we have $-\Gamma + ND^{-1}N\trans\in\mptnl(G_1\cup G_2)$.  By \cref{prop:psd} and the fact that $G_1\cup G_2$ is connected, the matrix $\Gamma - ND^{-1}N\trans$ is negative semidefinite with $n_0 = 1$.  Along with the fact that $D$ is positive definite, we have $n_-(B) = |V(G_1)| - 1$ and $n_0(B) = 1$.
\end{proof}

% WEAKDEF

Next, we define the weak nodal domain.  In the classical discrete nodal domain theorem, strong domains naturally form two groups, one with positive entries and one with negative entries.  In contrast, the strong domains for a general matrix $A\in\mptn(G)$ and its eigenvector $\by$ seem indistinguishable simply by signs.  However, we emphasize two features of weak nodal domains:  
\begin{enumerate}
    \item The weak nodal domains form a partition of the strong domains.
    \item Each zero domain is incident to at most two weak nodal domains.
\end{enumerate}    
For $A\in\mdso(G)$, every zero domain is adjacent to exactly two weak nodal domains, one positive and one negative.  However, \cref{ex:wdone} shows that this might not work for general matrices $A\in\mptn(G)$.

\begin{figure}[h]
\begin{center}
\begin{tikzpicture}
\begin{scope}[xshift=-3cm]
\node[label={below:$1$}] (1) at (0,0) {};
\foreach \i/\ang/\s in {2/70/+,3/110/-,4/190/+,5/230/-,6/310/+,7/350/-} {
    \node[label={\ang:$\i$}] (\i) at (\ang:1.5) {$\s$};
    \draw (1) --node[edgesign]{$+$} (\i);
}
\draw[very thick] (2) --node[edgesign]{$+$} (3);
\draw[very thick] (4) --node[edgesign]{$+$} (5);
\draw[very thick] (6) --node[edgesign]{$+$} (7);

\node[textonly] at (0,-1.5) {$G$};
\end{scope}

\begin{scope}[xshift=3cm]
\node[label={below:$z_1$}] (z1) at (0,0) {};
\node[label={90:$s_{23}$}] (s23) at (90:1.5) {};
\node[label={210:$s_{45}$}] (s45) at (210:1.5) {};
\node[label={330:$s_{67}$}] (s67) at (330:1.5) {};

\draw (z1) -- (s23);
\draw (z1) -- (s45);
\draw (z1) -- (s67);

\node[textonly] at (0,-1.5) {$G_{A,\by}$};
\end{scope}
\end{tikzpicture}
\end{center}
\caption{A graph $G$ and the corresponding $G_{A,\by}$.}
\label{fig:wdone}
\end{figure}
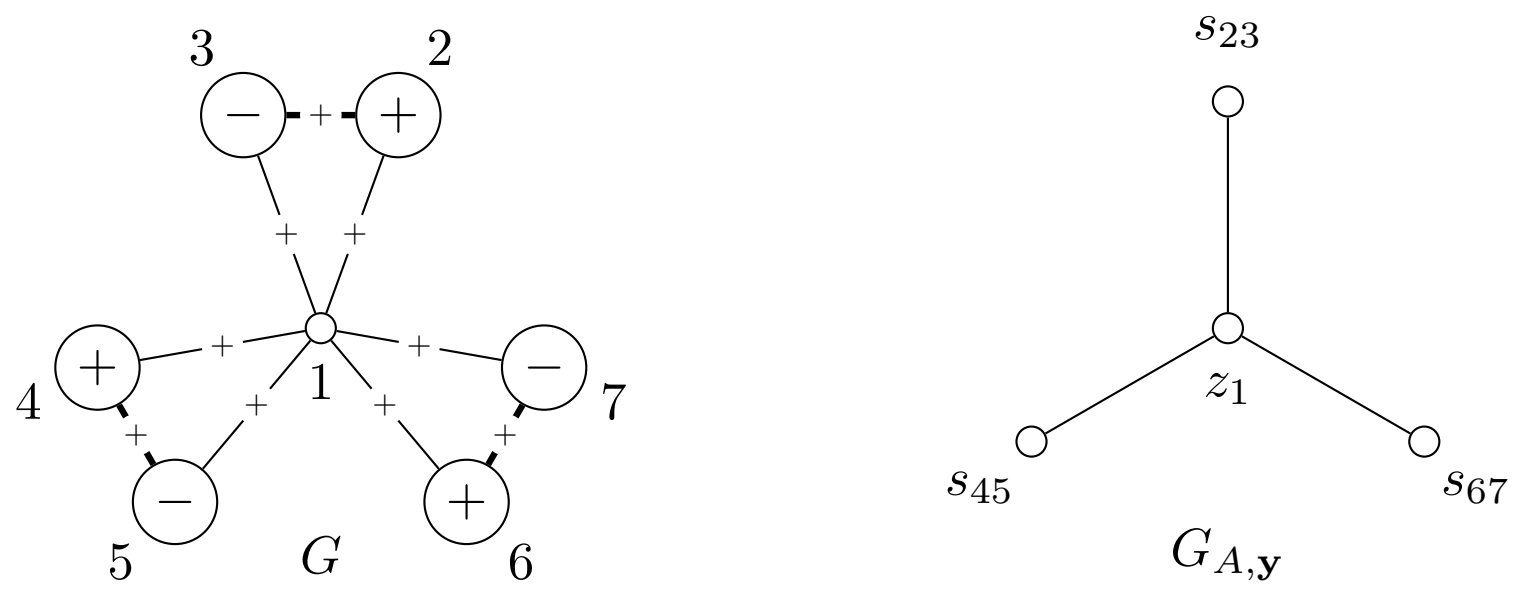

\begin{example}
\label{ex:wdone}
Let $G$ be the graph as shown in \cref{fig:wdone}.  Let $A$ be the matrix obtained from the adjacency matrix of $G$ by updating its $1,1$-entry to $6$.  Then $\spec(A) = \{-1, -1, -1, 0, 1, 1, 7\}$ and 
\[
    \by = (0,1,-1,1,-1,1,-1)\trans
\]
is an eigenvector of $A$ with respect to $\lambda = -1$.
In \cref{fig:wdone}, homogeneous edges are highlighted by thick edges, and the signs of $\by$ are marked inside the vertices.  Therefore, there are three strong domains, $\{2,3\}$, $\{4,5\}$, and $\{6,7\}$, and one zero domain, $\{1\}$.  However, with $\lambda = -1$, $n_-(A - \lambda I) + 1 = 1$ suggests that we need to define the weak nodal domains so that the zero domain $\{1\}$ merges all adjacent strong domains into one weak nodal domain.
\end{example}

\begin{definition}
Let $G$ be a graph, $A\in\mptn(G)$, and $\by\in\mathbb{R}^{V(G)}$.  The graph $G_{A,\by}$ is a minor of $G$ obtained by contracting all homogeneous edges and all edges with both endpoints in $V_{\by = 0}$.
\end{definition}

By definition, $G_{A - \lambda I,\by} = G_{A,\by}$.

\begin{remark}
Note that, for any given $A$ and $\by$, the strong domains and the zero domains form a partition of $V(G)$.  Naturally, each vertex of $G_{A,\by}$ represents a strong domain or a zero domain, and two vertices are adjacent if the corresponding two domains are adjacent by some edge.  
\end{remark}

\begin{definition}
The \emph{minor map} of $G_{A,\by}$ is a map $\phi: V(G) \rightarrow V(G_{A,\by})$ that sends $v$ to $i$, where $v$ is in the domain represented by $i$.  Let $V_*(G_{A,\by})$ be the set of vertices representing a strong domain, and $V_0(G_{A,\by})$ the set of vertices representing a zero domain.  
\end{definition}

In fact, a minor map $\phi: V(G) \rightarrow V(H)$ exists whenever $H$ is a minor of a connected graph $G$.  We refer the reader to \cite[Section~1.7]{Diestel25} for more details.

\begin{definition}
\label{def:spartition}
Let $G$ be a graph and $A\in\mptn(G)$.  Let $\by\in\mathbb{R}^{V(G)}$ and $\phi$ the minor map of $G_{A,\by}$.  For $z\in V_0(G_{A,\by})$, a \emph{signed partition} of $N_{G_{A,\by}}(z)$ is a partition $X\dunion Y = N_{G_{A,\by}}(z)$ such that there is a vertex $v\in\phi^{-1}(z)$ such that $s(v,\phi^{-1}(X))$ or $s(v,\phi^{-1}(Y))$ is nonzero.
\end{definition}

\begin{remark}
\label{rem:spartition}
Suppose $\by$ is an eigenvector of $A$ with respect to $\lambda$.  Let $z\in V_0(G_{A,\by})$, $v\in\phi^{-1}(z)$, and $X\dunion Y$ a partition of $N_{G_{A,\by}}(v)$.  Then $N_G(v)\subseteq \phi^{-1}(X\dunion Y\dunion \{z\})$.  
% Let $A_\lambda = A - \lambda I$.  
Since $\by$ is zero on $\phi^{-1}(z)$, the $v$-th entry of $(A - \lambda I) \by = \bzero$ can be written as
\[
    \sum_{u\in\phi^{-1}(X)} (A)_{v,u}(\by)_u + 
    \sum_{u\in\phi^{-1}(Y)} (A)_{v,u}(\by)_u
    = 0.
\]
If $X\dunion Y$ is a signed partition of $N_{G_{A,\by}}(z)$ witnessed by $s(v,X)$ being nonzero, then the first summation is nonzero, and necessarily the second summation is also nonzero since their sum is zero.
\end{remark}

\begin{definition}
\label{def:wdomain}
Let $G$ be a graph and $A\in\mptn(G)$.  Let $\by\in\mathbb{R}^{V(G)}$ and $\phi$ the minor map of $G_{A,\by}$.  A \emph{weak cover} $\mathcal{C}$ of $G_{A,\by}$ is a collection of its subgraphs of $G_{A,\by}$ that satisfies the following.  Here $\mathcal{C}_d = \{C\in\mathcal{C}: d\in V(C)\}$ for $d\in V(G_{A,\by})$.
\begin{enumerate}
\item For each $C\in\mathcal{C}$, $C$ is connected (including the case $C = K_1$) and every edge of $C$, if any, is incident to a vertex in $V_0(G_{A,\by})$. \label{wd-item:cconn}
\item For each $s\in V_*(G_{A,\by})$, $|\mathcal{C}_s| = 1$. \label{wd-item:s}
\item For each $z\in V_0(G_{A,\by})$, exactly one of the following holds. \label{wd-item:z}
\begin{enumerate}
    \item $\mathcal{C}_z = \emptyset$.
    \item $\mathcal{C}_z = \{C\}$ for some $C$, and $N_C(z) = N_{G_{A,\by}}(z)$.
    \item $\mathcal{C}_z = \{C_1, C_2\}$ for some $C_1,C_2$, and $N_{C_1}(z)\dunion N_{C_2}(z)$ is a signed partition of $N_{G_{A,\by}}(z)$.
\end{enumerate}
\item For each $C\in\mathcal{C}$, if $V(C)$ contains vertices $s\in V_*(G_{A,\by})$ and $z\in V_0(G_{A,\by})$ with $\{s,z\}\in E(G_{A,\by})$, then $\{s,z\}\in E(C)$. \label{wd-item:consistent}
\item The graph $G_{A,\by} - \{z: \mathcal{C}_z = \emptyset\}$ is connected. \label{wd-item:gconn}
\end{enumerate}
A \emph{weak domain combination} is a set $\{\phi^{-1}(V(C)): C\in\mathcal{C}\}$ for some weak cover $\mathcal{C}$, where an element in a weak domain combination is called a \emph{weak (nodal) domain}.  The maximum cardinality of a weak domain combination is denoted by $\nw(A,\by)$, or $\nw(\by)$ if $A$ is clear from the context.
\end{definition}

\begin{remark}
When $A\in\mdso(G)$ is a discrete Schr\"odinger operator, we claim that the weak nodal domains in the classical sense as in \cite{DGLS01} form a weak domain combination.  Recall that for a given $\by$, $V_+$, $V_-$, $V_0$ stands for the sets of vertices $v\in V(G)$ such that $(\by)_v$ is positive, negative, and zero, respectively.  Let $X_1, \ldots, X_q$ be the weak nodal domain in the classical sense, that is, they are the vertex sets of the components of $G[V_+\dunion V_0]$ or $G[V_+\dunion V_0]$.  For each $k = 1, \ldots, q$, one may define $C_k$ as the graph obtained from $G[X_k]$ by contracting all strong domains and zero domains and then removing edges between strong domains.  Thus, $\{C_k: k = 1, \ldots, q\}$ is a weak cover of $G_{A,\by}$.
\end{remark}

We will show $\nw(A,\by) \leq n_-(A - \lambda I) + 1$ in \cref{thm:jwndt}.  We first offer some examples for \cref{def:wdomain}.

\begin{example}
Let $G$, $A$, and $\by$ be the same as in \cref{ex:wdone}.  Then $G_{A,\by}$ is shown in \cref{fig:wdone}, and 
% \[
%     \begin{aligned}
%         V(G_{A,\by}) &= \{z_1, s_{23}, s_{45}, s_{67}\}, \\
%         E(G_{A,\by}) &= \{\{z_1,s_{23}\}, \{z_1,s_{45}\}, \{z_1,s_{67}\}\}, \\
%         V_0(G_{A,\by}) &= \{z_1\}, \text{ and} \\
%         V_*(G_{A,\by}) &= \{s_{23}, s_{45}, s_{67}\}.
%     \end{aligned}
% \]
\[
    V_0(G_{A,\by}) = \{z_1\} \text{ and } V_*(G_{A,\by}) = \{s_{23}, s_{45}, s_{67}\}.
\]
If $\mathcal{C}$ is a weak cover of $G_{A,\by}$, then $|\mathcal{C}_{z_1}|$ can be $0$, $1$, or $2$ by \cref{def:wdomain}\eqref{wd-item:z}.  If $|\mathcal{C}_{z_1}| = 0$, then the $G_{A,\by} - \{z: \mathcal{C}_z = \emptyset\}$ is not connected, violating \cref{def:wdomain}\eqref{wd-item:gconn}.  On the other hand, $N_{G_{A,\by}}(z_1)$ does not admit any signed partition since $s(1,\{2,3\}) = s(1,\{4,5\}) = s(1,\{6,7\}) = \circledast$.  Thus, $|\mathcal{C}_{z_1}| = 2$ is not an option.  Consequently, $|\mathcal{C}_{z_1}| = 1$.  This means $\mathcal{C} = \{G_{A,\by}\}$ is the only weak cover and $\nw(G,\by) = 1$.  Recall that $n_-(A - \lambda I,\by) + 1 = 1$ with $\lambda = -1$.
\end{example}

Next we provide an example where we exhaustively list all possible weak covers.  

\begin{figure}[h]
\begin{center}
\begin{tikzpicture}
\begin{scope}[xshift=-3cm]
\foreach \i/\ang in {1/90,2/210,3/330} {
    \node[label={\ang:$\i$}] (\i) at (\ang:2) {};
}
\node[label={150:$4$}] (4) at (150:2) {$+$};
\node[label={330:$5$}] (5) at (150:1) {$-$};
\node[label={270:$6$}] (6) at (270:2) {$+$};
\node[label={90:$7$}] (7) at (270:1) {$+$};
\node[label={30:$8$}] (8) at (30:2) {$+$};
\node[label={210:$9$}] (9) at (30:1) {$-$};

\draw (1) --node[edgesign]{$+$} (4) --node[edgesign]{$+$} (2) --node[edgesign]{$+$} (6) --node[edgesign]{$+$} (3) --node[edgesign]{$+$} (8) --node[edgesign]{$+$} (1);
\draw (1) --node[edgesign]{$+$} (5) --node[edgesign]{$+$} (2) --node[edgesign]{$-$} (7) --node[edgesign]{$-$} (3) --node[edgesign]{$+$} (9) --node[edgesign]{$+$} (1);
\draw[very thick] (4) --node[edgesign]{$+$} (5);
\draw (6) --node[edgesign]{$+$} (7);
\draw[very thick] (8) --node[edgesign]{$+$} (9);
\node[textonly] at (0,-3) {$G$};
\end{scope}

\begin{scope}[xshift=3cm]
\node[label={90:$z_1$}] (z1) at (90:2) {};
\node[label={210:$z_2$}] (z2) at (210:2) {};
\node[label={330:$z_3$}] (z3) at (330:2) {};
\node[label={330:$s_{45}$}] (s45) at (150:1) {};
\node[label={270:$s_6$}] (s6) at (270:2) {};
\node[label={90:$s_7$}] (s7) at (270:1) {};
\node[label={210:$s_{89}$}] (s89) at (30:1) {};

\draw (z1) -- (s45) -- (z2) -- (s7) -- (z3) -- (s89) -- (z1);
\draw (z2) -- (s6) -- (z3);
\draw (s6) -- (s7);

\node[textonly] at (0,-3) {$G_{A,\by}$};
\end{scope}
\end{tikzpicture}
\end{center}
\caption{A graph $G$ and the corresponding $G_{A,\by}$.}
\label{fig:allwc1}
\end{figure}
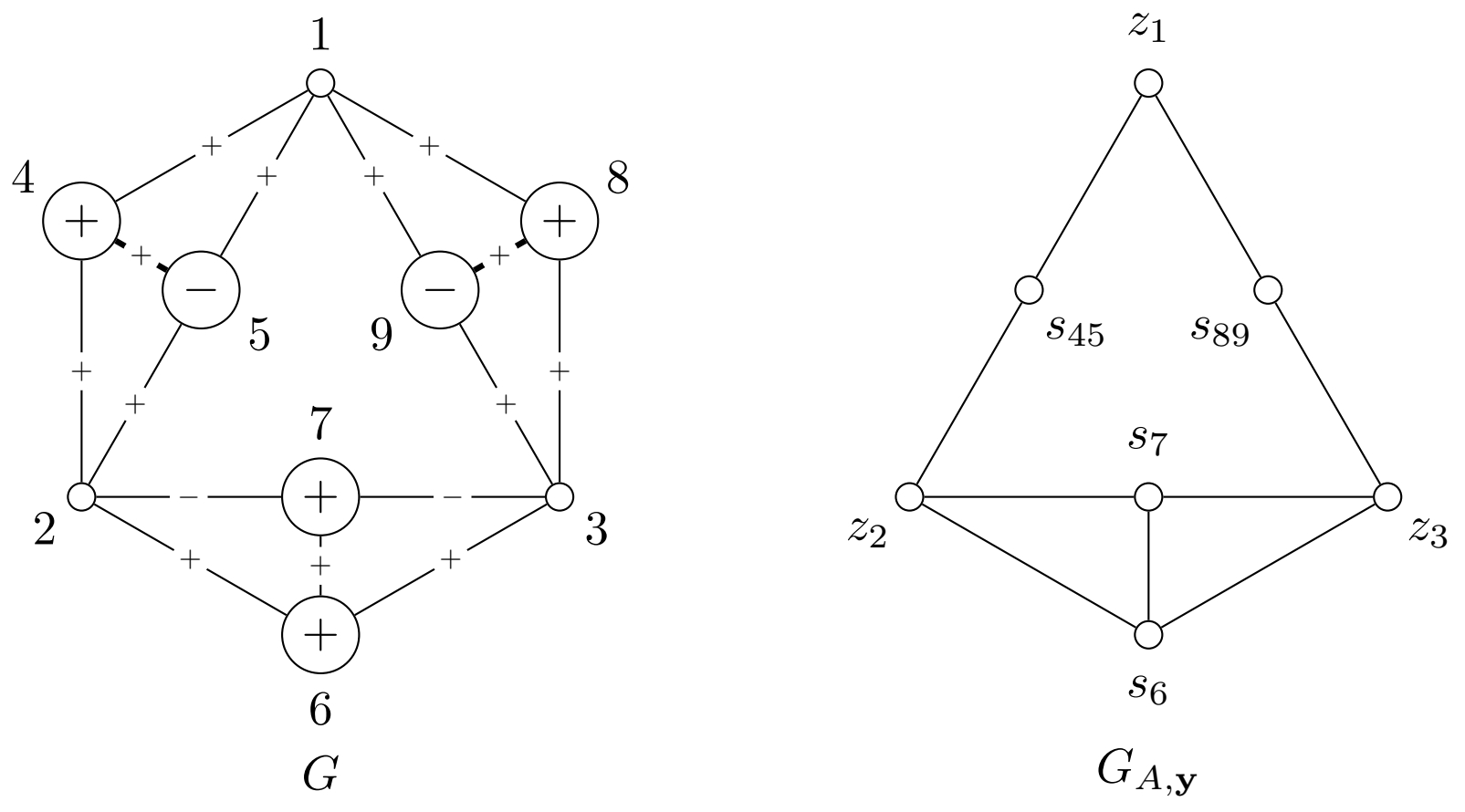

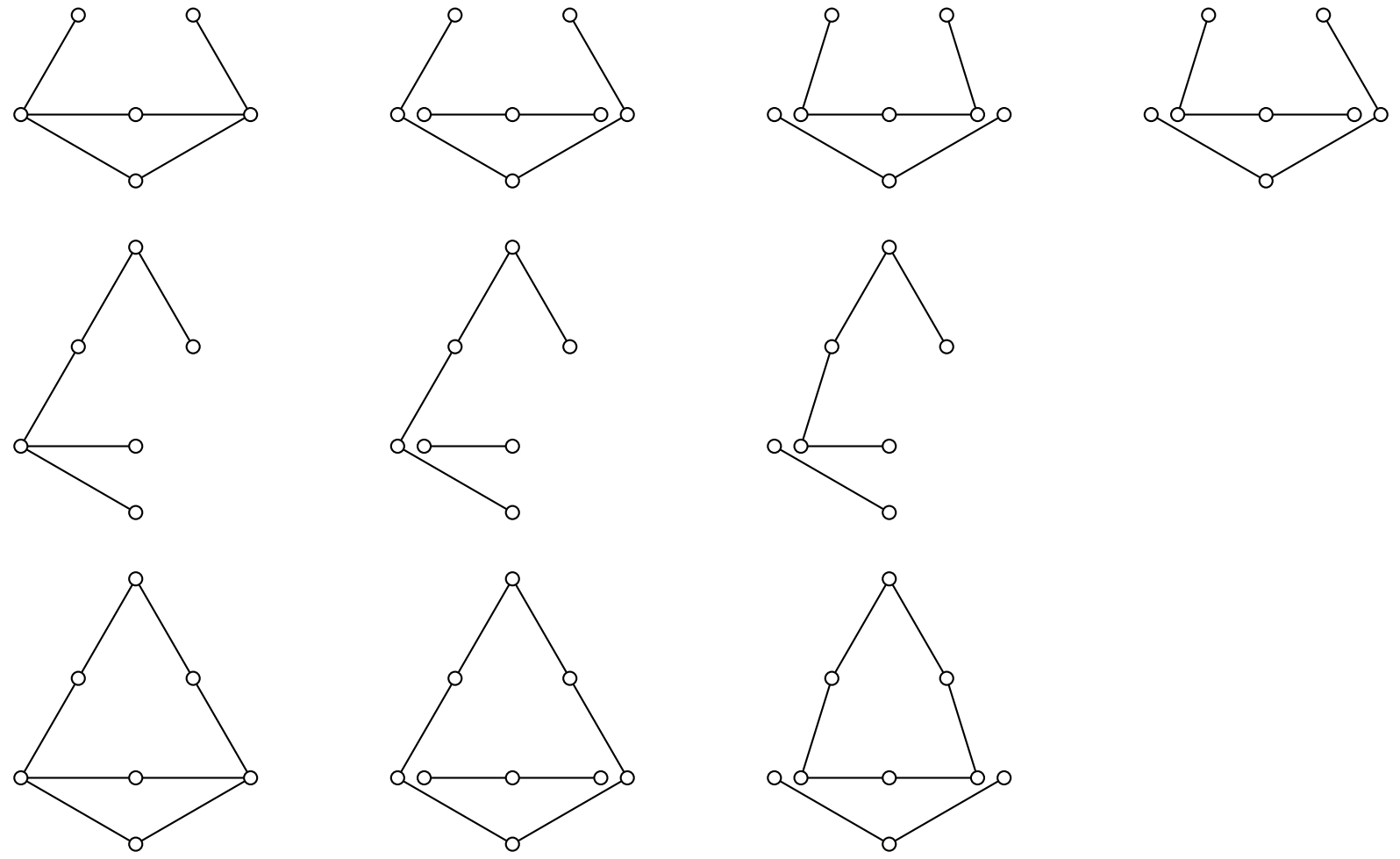
\begin{figure}[h]
%%%%%%%%%% row 1 %%%%%%%%%%
\begin{center}
% -z1: 1
\begin{tikzpicture}[transform shape, scale=0.5]
\node (z2) at (210:2) {};
\node (z3) at (330:2) {};
\node (s45) at (150:1) {};
\node (s6) at (270:2) {};
\node (s7) at (270:1) {};
\node (s89) at (30:1) {};

\draw (s45) -- (z2) -- (s7) -- (z3) -- (s89);
\draw (z2) -- (s6) -- (z3);
\end{tikzpicture}
\hfil
% -z1: 2
\begin{tikzpicture}[transform shape, scale=0.5]
\node (z2) at (210:2) {};
\node (z2p) at ([xshift=0.4cm] z2.center) {};
\node (z3) at (330:2) {};
\node (z3p) at ([xshift=-0.4cm] z3.center) {};
\node (s45) at (150:1) {};
\node (s6) at (270:2) {};
\node (s7) at (270:1) {};
\node (s89) at (30:1) {};

\draw(s45) -- (z2) -- (s6) -- (z3) -- (s89);
\draw (z2p) -- (s7) -- (z3p);
\end{tikzpicture}
\hfil
% -z1: 3
\begin{tikzpicture}[transform shape, scale=0.5]
\node (z2) at (210:2) {};
\node (z2p) at ([xshift=0.4cm] z2.center) {};
\node (z3) at (330:2) {};
\node (z3p) at ([xshift=-0.4cm] z3.center) {};
\node (s45) at (150:1) {};
\node (s6) at (270:2) {};
\node (s7) at (270:1) {};
\node (s89) at (30:1) {};

\draw(s45) -- (z2p) -- (s7) -- (z3p) -- (s89);
\draw (z2) -- (s6) -- (z3);
\end{tikzpicture}
\hfil
% -z1: 4
\begin{tikzpicture}[transform shape, scale=0.5]
\node (z2) at (210:2) {};
\node (z2p) at ([xshift=0.4cm] z2.center) {};
\node (z3) at (330:2) {};
\node (z3p) at ([xshift=-0.4cm] z3.center) {};
\node (s45) at (150:1) {};
\node (s6) at (270:2) {};
\node (s7) at (270:1) {};
\node (s89) at (30:1) {};

\draw (s45) -- (z2p) -- (s7) -- (z3p);
\draw (z2) -- (s6) -- (z3) -- (s89);
\end{tikzpicture}
\end{center}
%%%%%%%%%% row 2 %%%%%%%%%%
\begin{center}
% -z3: 1
\begin{tikzpicture}[transform shape, scale=0.5]
\node (z1) at (90:2) {};
\node (z2) at (210:2) {};
\node[opacity=0] (z3) at (330:2) {};
\node (s45) at (150:1) {};
\node (s6) at (270:2) {};
\node (s7) at (270:1) {};
\node (s89) at (30:1) {};

\draw (z1) -- (s45) -- (z2) -- (s7);
\draw (s89) -- (z1);
\draw (z2) -- (s6);
\end{tikzpicture}
\hfil
% -z3: 2
\begin{tikzpicture}[transform shape, scale=0.5]
\node (z1) at (90:2) {};
\node (z2) at (210:2) {};
\node (z2p) at ([xshift=0.4cm] z2.center) {};
\node[opacity=0] (z3) at (330:2) {};
\node (s45) at (150:1) {};
\node (s6) at (270:2) {};
\node (s7) at (270:1) {};
\node (s89) at (30:1) {};

\draw (z1) -- (s45) -- (z2) -- (s6);
\draw (s89) -- (z1);
\draw (z2p) -- (s7);
\end{tikzpicture}
\hfil
% -z3: 3
\begin{tikzpicture}[transform shape, scale=0.5]
\node (z1) at (90:2) {};
\node (z2) at (210:2) {};
\node (z2p) at ([xshift=0.4cm] z2.center) {};
\node[opacity=0] (z3) at (330:2) {};
\node (s45) at (150:1) {};
\node (s6) at (270:2) {};
\node (s7) at (270:1) {};
\node (s89) at (30:1) {};

\draw (z1) -- (s45) -- (z2p) -- (s7);
\draw (s89) -- (z1);
\draw (z2) -- (s6);
\end{tikzpicture}
\hfil
% gap
\begin{tikzpicture}[opacity=0, transform shape, scale=0.5]
\node (z1) at (90:2) {};
\node (z2) at (210:2) {};
\node[opacity=0] (z3) at (330:2) {};
\node (s45) at (150:1) {};
\node (s6) at (270:2) {};
\node (s7) at (270:1) {};
\node (s89) at (30:1) {};

\draw (z1) -- (s45) -- (z2) -- (s7);
\draw (s89) -- (z1);
\draw (z2) -- (s6);
\end{tikzpicture}
\end{center}
%%%%%%%%%% row 3 %%%%%%%%%%
\begin{center}
% -0:1
\begin{tikzpicture}[transform shape, scale=0.5]
\node (z1) at (90:2) {};
\node (z2) at (210:2) {};
\node (z3) at (330:2) {};
\node (s45) at (150:1) {};
\node (s6) at (270:2) {};
\node (s7) at (270:1) {};
\node (s89) at (30:1) {};

\draw (z1) -- (s45) -- (z2) -- (s7) -- (z3) -- (s89) -- (z1);
\draw (z2) -- (s6) -- (z3);
\end{tikzpicture}
\hfil
% -0:2
\begin{tikzpicture}[transform shape, scale=0.5]
\node (z1) at (90:2) {};
\node (z2) at (210:2) {};
\node (z2p) at ([xshift=0.4cm] z2.center) {};
\node (z3) at (330:2) {};
\node (z3p) at ([xshift=-0.4cm] z3.center) {};
\node (s45) at (150:1) {};
\node (s6) at (270:2) {};
\node (s7) at (270:1) {};
\node (s89) at (30:1) {};

\draw (z1) -- (s45) -- (z2) -- (s6) -- (z3) -- (s89) -- (z1);
\draw (z2p) -- (s7) -- (z3p);
\end{tikzpicture}
\hfil
% -0:3
\begin{tikzpicture}[transform shape, scale=0.5]
\node (z1) at (90:2) {};
\node (z2) at (210:2) {};
\node (z2p) at ([xshift=0.4cm] z2.center) {};
\node (z3) at (330:2) {};
\node (z3p) at ([xshift=-0.4cm] z3.center) {};
\node (s45) at (150:1) {};
\node (s6) at (270:2) {};
\node (s7) at (270:1) {};
\node (s89) at (30:1) {};

\draw (z1) -- (s45) -- (z2p) -- (s7) -- (z3p) -- (s89) -- (z1);
\draw (z2) -- (s6) -- (z3);
\end{tikzpicture}
\hfil
% gap
\begin{tikzpicture}[opacity=0, transform shape, scale=0.5]
\node (z1) at (90:2) {};
\node (z2) at (210:2) {};
\node (z3) at (330:2) {};
\node (s45) at (150:1) {};
\node (s6) at (270:2) {};
\node (s7) at (270:1) {};
\node (s89) at (30:1) {};

\draw (z1) -- (s45) -- (z2) -- (s7) -- (z3) -- (s89) -- (z1);
\draw (z2) -- (s6) -- (z3);
\draw (s6) -- (s7);
\end{tikzpicture}
\end{center}

\caption{For $G_{A,\by}$ in \cref{fig:allwc1}, there are $14$ weak covers.  Here we list $10$ of them, while the other $4$ can be obtained by left-right flipping the asymmetric cases.}
\label{fig:allwc2}
\end{figure}

\begin{example}
\label{ex:allwc}
Let $G$ be as shown in \cref{fig:allwc1}.  Let 
\[
    A = \begin{bmatrix}
  10 & 0 & 0 & 1 & 1 & 0 & 0 & 1 & 1 \\
  0 & 10 & 0 & 1 & 1 & 1 & -1 & 0 &  0 \\
  0 & 0 & 10 & 0 & 0 & 1 & -1 & 1 & 1 \\
  1 & 1 & 0 & 1 & 1 & 0 & 0 & 0 & 0 \\
  1 & 1 & 0 & 1 & 1 & 0 & 0 & 0 & 0 \\
  0 & 1 & 1 & 0 & 0 & -1 & 1 & 0 & 0 \\
  0 & -1 & -1 & 0 & 0 & 1 & -1 & 0 & 0 \\
  1 & 0 & 1 & 0 & 0 & 0 & 0 & 1 & 1 \\
  1 & 0 & 1 & 0 & 0 & 0 & 0 & 1 & 1
    \end{bmatrix}
    \text{ and }
    \by = \begin{bmatrix}
        0 \\ 0 \\ 0 \\ 1 \\ -1 \\ 1 \\ 1 \\ 1 \\ -1
    \end{bmatrix}.
\]
Then $\by$ is an eigenvector of $A$ with respect to $\lambda = 0$.  In \cref{fig:allwc1}, the signs of $A$ are marked on the edges, and the signs of $\by$ are marked inside the vertices.  There are two homogeneous edges, $\{4,5\}$ and $\{8,9\}$, so the strong domains are $\{4,5\}$, $\{6\}$, $\{7\}$, $\{8,9\}$ and the corresponding $G_{A,\by}$ is also shown in \cref{fig:allwc1}.  

Let $\phi$ be the minor map.  To find all weak covers of $G_{A,\by}$, we first examine the possibilities of $\mathcal{C}_{z_1}$, $\mathcal{C}_{z_2}$, and $\mathcal{C}_{z_3}$.  Observe that $|\mathcal{C}_{z_1}| = 0, 1$ since $N_{G_{A,\by}}(z_1)$ does not admit any signed partition.  For $\mathcal{C}_{z_2}$ and $\mathcal{C}_{z_3}$, its cardinality can be $0$, $1$, or $2$.  Also, $N_{G_{A,\by}}(z_2)$ admits two signed partitions, $\{s_{45},s_6\}\cup\{s_7\}$ and $\{s_6\}\cup\{s_{45},s_7\}$.  Similarly, $N_{G_{A,\by}}(z_3)$ admits signed partitions $\{s_{89},s_6\}\cup\{s_7\}$ and $\{s_6\}\cup\{s_{89},s_7\}$.

Next, we see that $\{z: \mathcal{C}_z = \emptyset\}$ contains at most one of $z_1,z_2,z_3$, for otherwise $G_{A,\by} - \{z: \mathcal{C}_z = \emptyset\}$ is disconnected.  In \cref{fig:allwc2}, the weak covers with $\{z: \mathcal{C}_z = \emptyset\} = \{z_1\}$ are listed in the first row, those with $\{z: \mathcal{C}_z = \emptyset\} = \{z_3\}$ are listed in the second row, and those with $\{z: \mathcal{C}_z = \emptyset\} = \emptyset$ are listed in the third row.  Note that the weak covers with $\{z: \mathcal{C}_z = \emptyset\} = \{z_2\}$ are omitted since they can be obtained by symmetry.

Thus, we have obtained $\nw(A,\by) = 2$, and one may calculate $n_-(A - \lambda I) + 1 = 2$.
\end{example}

Mohammadian~\cite{Mohammadian16} defined the weak domain by setting $\mathcal{C}_z = \emptyset$ for every $z\in V_0(G_{A,\by})$.  Ge and Liu~\cite{GL23} extended this idea and gave another definition of weak domain.  We note that both of their definitions are special cases of our weak domains defined through weak covers.

% WEAKTHM

\begin{theorem}
\label{thm:jwndt}
Let $G$ be a connected graph and $A\in\mptn(G)$.  Let $\by$ be an eigenvector of $A$ with respect to $\lambda$.  Then 
\[
    \nw(A,\by) \leq n_- + 1,
\]
where $n_- = n_-(A - \lambda I)$.
\end{theorem}
\begin{proof}
Let $\phi$ be the minor map of $G_{A,\by}$.  Let $\mathcal{C} = \{C_1, \ldots, C_q\}$ be a weak cover with $q = \nw(A,\by)$.  Let $X_i = \phi^{-1}(V(C_i))$ be the weak domains for $i = 1, \ldots, q$.  Define $\by_i$ as the vector obtained from $\by$ by keeping the same value for entries in $X_i$ and setting it to zero for others.  

Let $z_1, \ldots, z_r$ be the vertices $z\in V_0(G_{A,\by})$ with $|\mathcal{C}_z| = 2$ and let $Z_i = \phi^{-1}(z_i)$ for $i = 1, \ldots, r$.  For each $i$, if the signed partition of $N_{G_{A,\by}}(z_i)$ is witnessed by $s(v,\phi^{-1}(X))$ for some $v\in Z_i$ and some $X\subseteq N_{G_{A,\by}}(z_i)$, then set $\bz_i$ as the vector in $\mathbb{R}^{V(G)}$ whose $v$-th entry is $1$ and zero otherwise.

Let $A_\lambda = A - \lambda I$.  Let $I_0\in\mathbb{R}^{V(G)\times V(G)}$ be the diagonal matrix whose $v,v$-entry is $1$ if $v\in V_{\by = 0}$ and $0$ otherwise.  Choose a large number $t > 0$ such that $(A_\lambda + tI_0)[V_{\by = 0}]$ is positive definite.
Let $S$ be the $n\times (r + q)$ matrix whose columns are $\bz_1, \ldots, \bz_r, \by_1, \ldots, \by_q$.  In the following, we claim $B = S(A_\lambda + tI_0)S\trans$ has the form 
\[
    B = \begin{bmatrix}
        D & N\trans \\
        N & \Gamma
    \end{bmatrix},
\]
where $D$ is an $r\times r$ diagonal matrix with positive diagonal entries, $N$ is a weighted incidence matrix of $G_2$, and $-\Gamma\in\mptnl(G_1)$.  Here $G_1$ and $G_2$ are defined by
\[
    \begin{aligned}
    V(G_1) &= V(G_2) = \{x_i: i = 1, \ldots, q\}, \\
    E(G_1) &= \{\{x_i,x_j\}: \text{ there is an edge between $X_i$ and $X_j$ in $G_*$}\}, \text{ and}\\
    E(G_2) &= \{\{x_i,x_j\}: \mathcal{C}_{z_k} = \{C_i, C_j\} \text{ for some }k = 1, \ldots, r\}.
    \end{aligned}
\]

Note that there are no edges between distinct zero domains, for otherwise they will merge into one zero domain, so $D$ is a diagonal matrix.  By our choice of $t$, $(A + tI_0)[V_{\by = 0}]$ is positive definite, so the diagonal entries of $D$ are positive.

Next, we focus on $N$ and $\Gamma$.  Observe that
\[
    \begin{bmatrix}
        D & N\trans \\
        N & \Gamma
    \end{bmatrix}
    \begin{bmatrix}
        \bzero \\ \bone
    \end{bmatrix}
    =
    S\trans (A_\lambda + tI_0) S
    \begin{bmatrix}
        \bzero \\ \bone
    \end{bmatrix}
    = S\trans (A_\lambda + tI_0) \by = \bzero,
\]
where the last equality holds because $I_0$ contains nonzero diagonal entries only on $V_{\by = 0}$.  This means $N\trans\bone = \bzero$ and $\Gamma\bone = \bzero$.  For each $k = 1, \ldots, r$, we know that $\mathcal{C}_{z_k} = \{C_i, C_j\}$ for some $i,j$, so the $k$-th column of $N$ contains at most two potentially nonzero entries, the $i$-th entry and the $j$-th entry.  If the signed partition of $N_{G_{A,\by}}(z_k)$ is witnessed by $s(v,\phi^{-1}(N_{C_i}(z_k)))$ being nonzero, then the $i,k$-entry of $N$ is nonzero, and the $j,k$-entry of $N$ is also nonzero by \cref{rem:spartition}.  Consequently, $N$ is a weighted incidence matrix of $G_2$.  On the other hand, $-\Gamma$ is a weighted Laplacian matrix of $G_1$ by the same arguments as in the proof of \cref{thm:jsndt}.  Additionally, $G_1\cup G_2$ is connected by \cref{def:wdomain}\eqref{wd-item:gconn}.

By \cref{lem:dnngamma}, $n_-(B) = q - 1$ and $n_0(B) = 1$.  By \cref{lem:congruent,lem:weyl}, 
\[
    q - 1 = n_-(B) \leq n_-(A_\lambda + tI_0) \leq n_-(A_\lambda) = n_-.
\]
Therefore, the desired inequality $\nw(A,\by) \leq n_- - 1$ holds.
\end{proof}

We end this section with an example, showing that our new definition of weak domain is flexible enough to provides a tighter lower bound for $n_-(A - \lambda I)$, even in the case of $A\in\mdso(G)$.

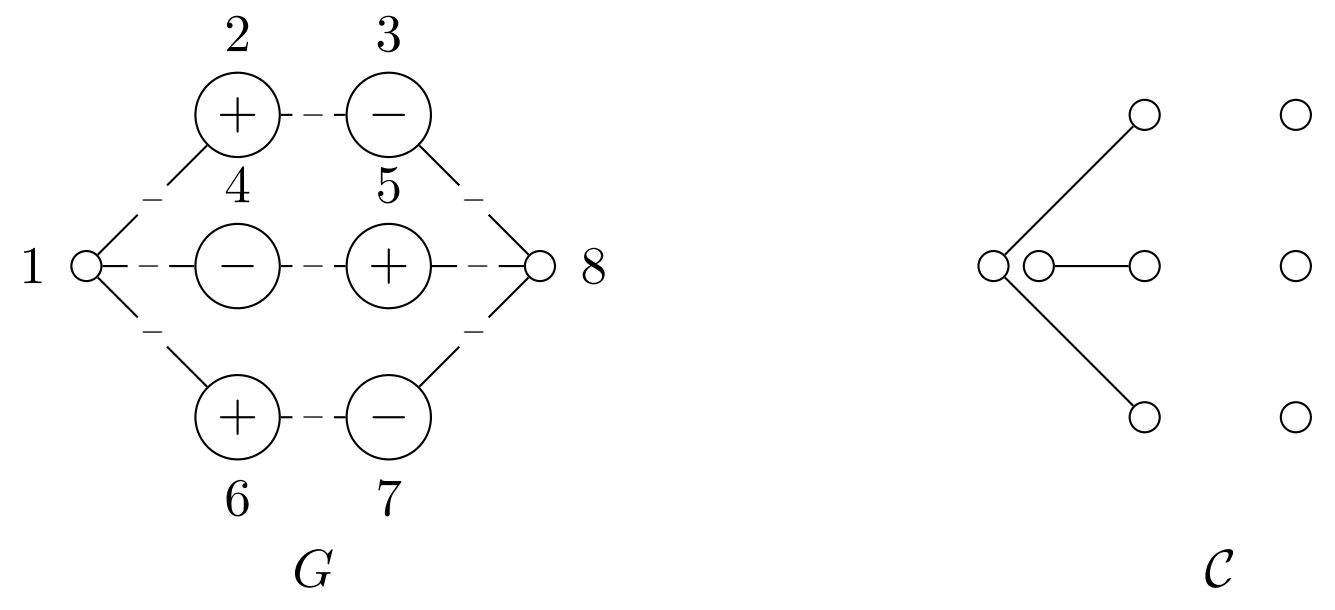
\begin{figure}[h]
\begin{center}
\begin{tikzpicture}
\begin{scope}[xshift=-3cm]
\node[label={left:$1$}] (1) at (0,0) {};
\node[label={above:$2$}] (2) at (1,1) {$+$};
\node[label={above:$3$}] (3) at (2,1) {$-$};
\node[label={above:$4$}] (4) at (1,0) {$-$};
\node[label={above:$5$}] (5) at (2,0) {$+$};
\node[label={below:$6$}] (6) at (1,-1) {$+$};
\node[label={below:$7$}] (7) at (2,-1) {$-$};
\node[label={right:$8$}] (8) at (3,0) {};
\draw (1) --node[edgesign]{$-$} (2) --node[edgesign]{$-$} (3) --node[edgesign]{$-$} (8);
\draw (1) --node[edgesign]{$-$} (4) --node[edgesign]{$-$} (5) --node[edgesign]{$-$} (8);
\draw (1) --node[edgesign]{$-$} (6) --node[edgesign]{$-$} (7) --node[edgesign]{$-$} (8);
\node[textonly] at (1.5,-2) {$G$};
\end{scope}[xshift=-3cm]

\begin{scope}[xshift=3cm]
\node (z1) at (0,0) {};
\node (z1p) at ([xshift=0.3cm]z1.center) {};
\node (s2) at (1,1) {};
\node (s3) at (2,1) {};
\node (s4) at (1,0) {};
\node (s5) at (2,0) {};
\node (s6) at (1,-1) {};
\node (s7) at (2,-1) {};
% \node (z8) at (3,0) {};
\draw (z1) -- (s2);
\draw (z1p) -- (s4);
\draw (z1) -- (s6);
\node[textonly] at (1.5,-2) {$\mathcal{C}$};
\end{scope}
\end{tikzpicture}
\end{center}
\caption{A graph $G$ and a weak cover $\mathcal{C}$.}
\label{fig:wdtight}
\end{figure}

\begin{example}
\label{ex:wdtight}
Let $G$ be as shown in \cref{fig:wdtight}.  Let 
\[
    A = \begin{bmatrix}
    0 & -1 & 0 & -1 & 0 & -1 & 0 & 0 \\
    -1 & 0 & -1 & 0 & 0 & 0 & 0 & 0 \\
    0 & -1 & 0 & 0 & 0 & 0 & 0 & -1 \\
    -1 & 0 & 0 & 0 & -1 & 0 & 0 & 0 \\
    0 & 0 & 0 & -1 & 0 & 0 & 0 & -1 \\
    -1 & 0 & 0 & 0 & 0 & 0 & -1 & 0 \\
    0 & 0 & 0 & 0 & 0 & -1 & 0 & -1 \\
    0 & 0 & -1 & 0 & -1 & 0 & -1 & 0
    \end{bmatrix}
    \text{ and }
    \by = \begin{bmatrix}
    0 \\ 1 \\ -1 \\ -2 \\ 2 \\ 1 \\ -1 \\ 0
    \end{bmatrix}.
\]
Then $\by$ is an eigenvector $A$ with respect to $\lambda = 1$, and $n_-(A - \lambda I) + 1 = 5$.  In view of the classical nodal domain theorem, the graph has $4$ weak nodal domains, so the bound is not tight.  However, the weak cover $\mathcal{C}$ in \cref{fig:wdtight} leads to a weak domain combination with $5$ weak domains, showing $\nw(A,\by) \geq 5$.  So \cref{thm:jwndt} offers a tight bound in this example.
\end{example}

% \section{NP-hardness}
% \label{sec:hard}

\section{Concluding remarks}
\label{sec:conclusion}

In this paper, we demonstrate that the congruent quotient matrix and its inequalities are the essence of the proof of the discrete nodal domain theorem.  The techniques allow us to extend the existing bounds for the number of strong domains and refine the definition of weak domain.  Consequently, we also enhance the bounds for $|\het(A,\by)|$ and obtain Fiedler's result (\cref{thm:f75a2}) as a corollary, connecting all related theorems.

While the new definition of weak domain seems artificial, it maintains the same bound $\nw(A,\by) \leq n_- + 1$ and generalizes the weak nodal domains in the classical sense as well as those defined in \cite{Mohammadian16, GL23}.  Even for $A\in\mdso(G)$, the new definition leads to an improvement to the bound of $n_-$.  

When the signs of $A$ and $\by$ are given, finding $\nw(A,\by)$ is a combinatorial optimization problem, and it seems to be NP-hard since iterating all possible weak covers is tedious.  With the bound $\nw(A,\by) \leq n_- + 1$, if a weak cover $\mathcal{C}$ and a realization of $A$ and $\by$ are found so that $|\mathcal{C}| = n_- + 1$, then the optimal value is found.

These new results offer a clearer picture of the connections between the inertia and the signs of an eigenvector, while many questions remain open.  Given the signs of $A$ and $\by$, are there realizations of $A$ and $\by$ such that the bounds are tight?  Is there an efficient algorithm for computing $\nw(A,\by)$?  Also, it would be nice to see if any lower bound for $\nw(A,\by)$ exists.

% \bibliography{./JLaTeX/AuthorA,./JLaTeX/JournalA,./JLaTeX/JepBib}{}
% \bibliographystyle{plain}

%bibitem code

\end{document}